\documentclass[10pt,reqno]{article}
\usepackage{longtable,latexsym,epsfig,caption,subcaption,microtype,textcomp,stmaryrd,listings,lineno,hyperref}

\usepackage{amssymb}
\usepackage{graphicx}
\usepackage{amscd}
\usepackage{amsthm}

\usepackage{float}
\usepackage{graphics, amsmath}
\usepackage{amsfonts}
\usepackage{latexsym}
\usepackage{epsf}
\usepackage{xcolor, color}
\begin{document}

\theoremstyle{plain}
\newtheorem{theorem}{Theorem}
\newtheorem{corollary}[theorem]{Corollary}
\newtheorem{lemma}{Lemma}
\newtheorem{proposition}[theorem]{Proposition}

\theoremstyle{definition}
\newtheorem{definition}[theorem]{Definition}
\newtheorem{example}[theorem]{Example}
\newtheorem{conjecture}[theorem]{Conjecture}

\theoremstyle{remark}
\newtheorem{remark}{Remark}
\begin{center}
{\Large{{\bf Normalized logistic wavelets: Applications to COVID-19 data in Italy }}}
\medskip
\end{center}
\leftline{}
\centerline{\bf Grzegorz Rz\c{a}dkowski }

\leftline{Department of Finance and Risk Management, Warsaw University of Technology, Narbutta 85, 02-524 Warsaw, Poland}
\leftline{e-mail: grzegorz.rzadkowski@pw.edu.pl}

\newcommand{\Eulerian}[2]{\genfrac{<}{>}{0pt}{}{#1}{#2}}

\begin{abstract}
In this paper we deal with the logistic wavelets introduced in \cite{RF}. We modify them by multiplying by appropriate coefficients so that their norm in the space $L^{2}(R)$ is equal to 1. We calculate the normalization coefficients using the Grosset-Veselov formula \cite{GV}, Eulerian numbers and Bernoulli numbers. Then we apply the logistic wavelets to model of the first wave of Covid-19 deaths in Italy in 2020. This example shows that even asymmetric and skewed data can be modeled, with high accuracy, by a sum of logistic functions.
\end{abstract}

Keywords: Logistic wavelet, logistic equation, logistic function, COVID-19, Eulerian number, Bernoulli number, Riccati's differential equation.

2020 Mathematics Subject Classification: 92D30, 65T60, 11B83

\section{Introduction} 
The logistic equation defining the logistic function  $x = x(t)$  has the form (cf. \cite{RF})
\begin{equation}\label{1a}
	x'(t)=\frac{s}{x_{max}}\:x(x_{max}-x),\quad x(0)=x_{0}.
\end{equation}
where $t$ is time, and parameters $s$-steepness or slope coefficient and $x_{max}$-saturation level are constants.
The integral curve $x(t)$ of equation  (\ref{1a}) satisfying the condition $0<x(t)<x_{max}$ is called the logistic function. The logistic function
is used to describe and model various phenomena in physics, economics, medicine, biology, engineering, sociology and many
other sciences. Logistic functions now seem even more important from the point of view of their possible applications, due to the theory of the Triple Helix (TH) developed in the 1990s by Etzkowitz and Leydesdorff \cite{EL} (see also Leydesdorff \cite{L}).  This theory explains the phenomenon of creating and introducing innovations under the influence of the interaction of three factors University-Industry-Government and relations between them. According to the TH theory, the phenomenon of the emergence of innovations can be described by means of  logistic functions. Ivanova \cite{Iv1}, \cite{Iv2}, \cite{Iv3} has shown that the KdV equation naturally appears in TH theory and has also applied it to other fields such as the COVID-19 pandemic or financial markets.

After solving the differential equation   (\ref{1a}) we obtain the logistic function in the form
\begin{equation}\label{1b}
	x(t)=\frac{x_{max}}{1+e^{-s(t-t_0)}},
\end{equation}
where $t_0$ is the inflection point associated with the initial condition $\displaystyle u(0)=u_{0}=\frac{u_{max}}{1+e^{st_0}}$, then $\displaystyle t_0=\frac{1}{s}\log\Big(\frac{x_{max}-x_{0}}{x_{0}}\Big)$. At the point $t_0$, $x(t_0)=x_{max}/2$.
Equation (\ref{1a})  is a special case of the Riccati equation with constant coefficients
\begin{equation}\label{1c}
	x'(t)=r(x-x_{1})(x-x_{2}),
\end{equation}
where constants $r\neq 0,\;x_{1},\;x_{2}$ can be real or more generally complex numbers. 

If $x=x(t)$ is the solution of  (\ref{1c}) then its $n$ derivative $x^{(n)}(t)$ ($n=2,3,4,\ldots$) is a polynomial of the function $x(t)$ \cite{Rz1}, \cite{Rz2}, \cite{F}
\begin{equation}\label{1d}
x^{(n)}(t) = r^{n}\sum\limits_{k=0}^{n-1}\Eulerian{ n}{k }
(x-x_{1})^{k+1}(x-x_{2})^{n-k}
\end{equation}
 for $n=2,3,\ldots $, where $\displaystyle \Eulerian{ n}{k }$ denotes Eulerian number (the number of permutations $\{1,2,\ldots ,n\}$ having exactly $k,\: (k=0,1,2,\ldots ,n-1)$ ascents, Graham et al \cite{GKP}. 

Formula (\ref{1d}) applied to the logistic equation  (\ref{1a}) yields:
\begin{equation}\label{2a}
	x^{(n)}(t)= \left(-\frac{s}{x_{max}} \right)^{n}\;\sum\limits_{k=0}^{n-1}\Eulerian{n}{ k}
x^{k+1}(x-x_{max})^{n-k}.
\end{equation}

The paper has the following structure. In Section~\ref{sec2} we first briefly describe the general wavelet theory and then the logistic wavelets introduced in article A. Then we compute the normalizing coefficients for them. Section~\ref{sec3} is devoted to an application of logistic wavelets to model the spread of the COVID-19 pandemic in Italy in 2020. The paper is concluded in Section~\ref{sec4}. All data used in the paper were obtained from the website Our World in Data \cite{Our}.

\section{Wavelets and logistic wavelets}\label{sec2}
\subsection{Wavelets}
Let us now recall some general facts about wavelet theory  (cf. \cite{D, MR, M}) which we will use later. A wavelet or mother wavelet
( Daubechies \cite{D}, p.24 ) is an integrable function  $\psi \in L^{1}(\mathbb{R})$ with the following admissibility condition:
\begin{equation}\label{3a}
	C_{\psi}=2\pi\int_{-\infty}^{\infty}|\xi|^{-1}|\widehat{\psi}(\xi)|^2 d\xi < \infty,
\end{equation}
where $\widehat{\psi}(\xi)$ is the Fourier transform of   $\psi$
\[
	\widehat{\psi}(\xi)=\frac{1}{\sqrt{2\pi}}\int_{-\infty}^{\infty}\psi(x)e^{-i\xi x} dx. \]
	
Since the function $\psi \in L^{1}(\mathbb{R})$, then $\widehat{\psi}(\xi)$ is a continuous function, and condition (\ref{3a}) is satisfied only when $\widehat{\psi}(0)=0$ or $\int_{-\infty}^{\infty}\psi(x)dx=0$. On the other hand, Daubechies \cite{D}, p.24  shows that condition $\int_{-\infty}^{\infty}\psi(x)dx=0$ together with the second condition, slightly stronger than integrability, namely  $\int_{-\infty}^{\infty}|\psi(x)|(1+|x|)^{\alpha}dx<\infty$, for some  $\alpha>0$ are sufficient for (\ref{3a}). Usually much more is assumed about the function  $\psi$,so from a practical point of view the conditions $\int_{-\infty}^{\infty}\psi(x)dx=0$ and (\ref{3a}) are equivalent.Suppose furthermore that $\psi$ is also square integrable, $\psi \in L^{2}(\mathbb{R})$ with the norm
\[ || \psi ||=\left(\int_{-\infty}^{\infty}|\psi(x)|^2 dx \right)^{1/2}. \]
Using the mother wavelet, by dilating and translating, a double-indexed family of wavelets is obtained
	\[\psi^{a,b}(x)=\frac{1}{\sqrt{|a|}}\psi\Big(\frac{x-b}{a}\Big),
\]
where $a,b\in \mathbb{R}, \; a\neq 0$. The normalization has been chosen so that $||\psi^{a,b}||=||\psi||$ for all $a, b$. In order to be able to compare different wavelet families with each other, it is usually assumed that $||\psi||=1 $. Continuous Wavelet Transform (CWT) of a function $f \in L^{2}(\mathbb{R})$ with respect to a given wavelet family is defined as
\begin{equation}\label{3bb}
(T^{wav}f)(a,b)=\langle f, \psi^{a,b} \rangle =\int_{-\infty}^{\infty}f(x) \psi^{a,b}(x) dx.
\end{equation}

\subsection{Logistic wavelets}\label{sec3}
Logistic mother wavelets understood as derivatives of the logistic function $x(t)=\frac{1}{1+e^{-t}}$, which is a solution to the logistic equation
\begin{equation}\label{3b}
x'(x)=x(1-x)=-x(x-1), 
\end{equation}
are described in \cite{RF}. We will use them here as well, but multiplied by appropriate factor so that their norms in the space $L^{2}(\mathbb{R})$ are equal to 1. Now we will show how these factors can be calculated. Obviously wavelets modified in this way also satisfy the admissibility condition (\ref{3a}). Formulas (\ref{1d}) or (\ref{2a}) applied to equation (\ref{3b}) give:

\begin{equation}\label{3c}
	x^{(n)}(t)=(-1)^{n}\sum\limits_{k=0}^{n-1}\Eulerian{ n}{k }x^{k+1}(x-1)^{n-k}=\sum\limits_{k=0}^{n-1}(-1)^{k}\Eulerian{ n}{k }x^{k+1}(1-x)^{n-k},
\end{equation}
for $n=2,3,\ldots$. We will use the following Grosset and Veselov formula \cite{GV}
\begin{equation}\label{3d}
	\int_{-\infty}^{+\infty} \left(\frac{d^{n-1}}{dt^{n-1}}
	\frac{1}{\cosh^{2}t}\right)^{2}dt=(-1)^{n-1}2^{2n+1} B_{2n},
\end{equation}
where $n=1,2,\ldots$ and $B_{2n}$ is the $2n$th Bernoulli number. Other proofs of the Grosset-Veselow formula can be found in \cite{Bo}, \cite{Rz3}. Bernoulli numbers have the following generating function  (see Graham, Knuth, Patashnik \cite{GKP})
	\[B(\xi) = B_0+B_{1}\xi + B_{2}\frac{\xi^{2}}{2!}+\cdots =\frac{\xi}{e^{\xi}-1},
	\qquad |\xi|<2\pi.\]
It is known that $B_{n}$ is zero for all odd numbers $n\ge 3$. These numbers are rational and occur in formulas such as
	\[\sum_{k=1}^{\infty}\frac{1}{k^{2n}}=(-1)^{n+1}
	\frac{2^{2n-1}\pi^{2n}}{(2n)!}B_{2n} \quad n=1,2,\ldots
\]
The first few Bernoulli numbers are as follows
	\[B_{0}=1,\;B_{1}=-\frac{1}{2},\;B_{2}=\frac{1}{6},\;B_{4}=-\frac{1}{30}
,\;B_{6}=\frac{1}{42},\;B_{8}=-\frac{1}{30},\;B_{10}=\frac{5}{66},\;B_{12}=-\frac{691}{2730}.
\]
Returning to the problem of normalizing the derivatives $x^{(n)}(t)$ of the function $x(t)=\frac{1}{1+e^{-t}}$, note that the integral (\ref{3d}) can be written in the following form ( we put $\tau=2t$ at the end)
\begin{align}
\int_{-\infty}^{+\infty} &\left(\frac{d^{n-1}}{dt^{n-1}}\frac{1}{\cosh^{2}t}\right)^{2}dt=
\int_{-\infty}^{+\infty} \left(\frac{d^{n-1}}{dt^{n-1}}\frac{4e^{-2t}}{(1+e^{-2t})^2}\right)^{2}dt=
4\int_{-\infty}^{+\infty} \left(\frac{d^{n}}{dt^{n}}\frac{1}{1+e^{-2t}}\right)^{2}dt \nonumber\\
&=4(2^n)^2\int_{-\infty}^{+\infty}(x^{(n)}(2t))^2dt=2(2^n)^2\int_{-\infty}^{+\infty}(x^{(n)}(\tau))^2d\tau.\label{3e}
\end{align}
Comparing (\ref{3e}) with  (\ref{3d}) we get
\begin{equation}\label{3f}
\int_{-\infty}^{+\infty}(x^{(n)}(t))^2dt=(-1)^{n-1}B_{2n}=|B_{2n}|.
\end{equation}
Then, based on (\ref{3f}), we can redefine, with respect to \cite{RF}, the logistic mother wavelet $\psi_n(t)$ of order $n=2, 3,\ldots$ as
\begin{equation}\label{3g}
\psi_n(t)=\frac{1}{\sqrt{|B_{2n}|}}x^{(n)}(t),
\end{equation}
with the norm $||\psi_n||=||\psi_n||_{L^2}=1$.

\begin{figure}
	\begin{center}
	 \includegraphics[height=6cm, width=8cm]{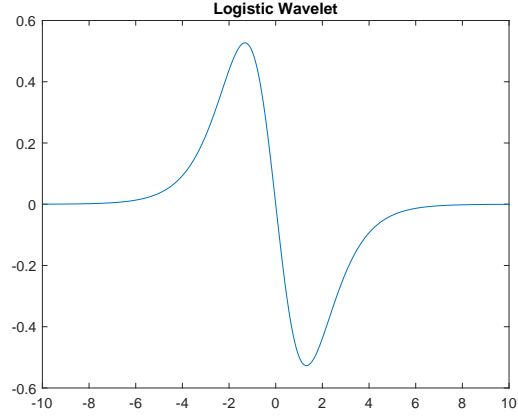}
	\end{center}
	\vspace{-7mm}
	\caption{Wavelet $\psi_2(t)$ }
	\label{fig1}
\end{figure}

In particular, for $n=2$ from the formula (\ref{3c}) or directly from (\ref{3b}) we have
\[x''(t)=x(1-x)(1-2x),
\]
and then by (\ref{3g}), wavelet $\psi_2(t)$  (Fig.~\ref{fig1}) is as follows
\begin{equation}\label{3h}
	\psi_2(t)=\frac{\sqrt{30}}{1+e^{-t}}\Big(1-\frac{1}{1+e^{-t}}\Big)\Big(1-\frac{2}{1+e^{-t}}\Big)=\frac{\sqrt{30}(e^{-2t}-e^{-t})}{(1+e^{-t})^3}.
\end{equation}

For $n=2$ or more generally for $n=2,3,\ldots$ we create,by dilating and translating, a doubly indexed family of wavelets (children wavelets)
	\[\psi_n^{a,b}(t)=\frac{1}{\sqrt{|a|}}\psi_n\Big(\frac{t-b}{a}\Big),
\]
where $a,b\in \mathbb{R}, \; a\neq 0$.

We implement the $\psi_2(t)$ wavelet in Matlab (Matlab's wavelet toolbox) with the following code:

function [psi,t] = logist(LB,UB,N,$\sim$)\\
\%\textit{LOGISTIC Logistic wavelet.}\\
\%   \textit{[PSI,T] = LOGIST(LB,UB,N) returns values of }\\
\%   \textit{the Logistic wavelet on an N point regular}\\
\%   \textit{grid in the interval [LB,UB].}\\
\%   \textit{Output arguments are the wavelet function PSI}\\
\%   \textit{computed on the grid T.}\\
\%   \textit{This wavelet has [-7 7] as effective support.}\\
\%   \textit{See also WAVEINFO.}\\
\% \textit{Compute values of the Logistic wavelet.}\\
t = linspace(LB,UB,N);        \% \textit{wavelet support.}\\
psi =sqrt(30)* (exp(-2*t)-exp(-t))./(1+exp(-t)).$^\wedge$3;\\
end

\section{An application to COVID-19 data in Italy}\label{sec3}
\begin{figure}
	\begin{center}
	 \includegraphics[height=6cm, width=8cm]{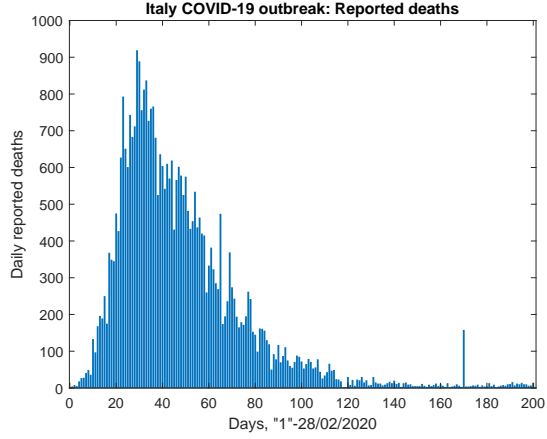}
	\end{center}
	\vspace{-7mm}
	\caption{The `first wave' of deaths in Italy, 28/02/2020--14/09/2020}
	\label{fig2}
\end{figure}

Let us consider the time series of daily deaths during the COVID-19 pandemic in Italy in the period from February 28, 2020 to September 14, 2020 (200 days), Fig.~\ref{fig2}. These data are known as the "first wave" of deaths in Italy and have already been analyzed many times by various authors (cf. e.g.,  Bezzini et al. \cite{B}, Dorrucci et al. \cite{DM}). It is clear that the time series Fig.~\ref{fig2} is not symmetric and is skewed to the right. It could be modeled by using, for example, the Gompertz function or another right-skewed distribution, which are broadly applied, e.g, in insurance \cite{HK}. We will show that, however, the time series of the total number of deaths can also be modeled with high accuracy by a sum of logistic functions. 

Let $(y_n)$ be the smoothed (by using 7-day moving averages) time series of the total, reported number of Covid-19 deaths in Italy, up to $n$th  day. Then we calculate its first differences, i.e., the daily numbers of deaths
\[\Delta^{1}y_n=y_n-y_{n-1},
\]
and the central second differences (changes in daily deaths)
\[\Delta^{2}y_n=\Delta^{1}y_{n+1}-\Delta^{1}y_n =y_{n+1}-  2 y_n+y_{n-1}.
\]

Assuming that $(y_n)$ follows locally a logistic function $\displaystyle y_n\approx y(n)=y_{max}/(1+\exp(-(n-b)/a))$ and applying definition (\ref{3h}) we have
\begin{equation}\label{4a}
	 y''(t)=\frac{y_{max}}{\sqrt{30}\cdot a^{3/2}}\psi_2^{a,b}(t).
\end{equation}
For the second differences $\Delta^{2}y_n$, we will apply the CWT transform (\ref{3bb}). Directly from the CWT scalogram we can find, for a given logistic wave, such values of the parameters $b$ and $a$ so that the value of the Index (\ref{4b}) at the point with coordinates $(b,a)$ is maximal. By (\ref{4a}) we have
\begin{align}
	\text{Index}= \sum\limits_{n}&\Delta^2y_n\psi_2^{a,b}(n)\approx \sum\limits_{n}\Delta^2y(n)\psi_2^{a,b}(n)
	\approx \int_{-\infty}^{\infty}y''(t)\psi_2^{a,b}(t)dt
	 =\int_{-\infty}^{\infty} \frac{y_{max}}{\sqrt{30}\cdot a^{3/2}}\psi_2^{a,b}(t)\psi_2^{a,b}(t)dt \nonumber\\
	&=\frac{y_{max}}{\sqrt{30}\cdot a^{3/2}}\int_{-\infty}^{\infty}(\psi_2^{a,b}(t))^2 dt= \frac{y_{max}}{\sqrt{30}\cdot a^{3/2}}.\label{4b}
\end{align}

Using (\ref{4b}) we can estimate the saturation level $y_{max}$ as follows
\begin{equation}\label{4c}
y_{max}\approx \sqrt{30}\cdot a^{3/2}\sum\limits_{n}\Delta^2y_n\psi_2^{a,b}(n)=\sqrt{30}\cdot a^{3/2}\text{Index}.
\end{equation}

Let us note that the parameter $b$  is the time (day) when a given logistic wave reached the inflection point (maximum of daily values). The parameter $a$ can be interpreted in terms of the length of a given logistic wave. Namely, for a logistic wave $x(t)$ of the form
\[x(t)=\frac{x_{max}}{1+\exp(-\frac{t-b}{a})}
\] 
we can define e.g., $95\%$ confidence interval by cutting off $2.5\%$ of the left and right values. Denoting by $t_1$ the left end of this interval we have
\[x(t_1)=\frac{x_{i,max}}{1+\exp(-\frac{t_1-b}{a})}=0.025x_{max},
\]
where, after an easy calculation, we get
\[b-t_1=3.66 a.\]
From the symmetry of the logistic function, the length of the confidence interval is $7.32a$. Thus, in practice, we can  assume that the length of this logistic wave is
\begin{equation}\label{4cc}
 \text{wavelength}=7.32a.
\end{equation}

Now we will model the time series $(y_n)$ by a sum of logistic functions
\begin{equation}\label{4bb}
f(t)=\sum_{i=1}^{k}\frac{x_{i,max}}{1+\exp(-\frac{t-b_i}{a_i})},
\end{equation}
$i=1,2,\ldots,k$, where $k$ is the number of logistic waves.

If there are several overlapping logistic waves, occuring in the same time period, then the higher intensity waves (with larger Index) may cause the lower intensity waves to be invisible on the CWT scalogram. Therefore, in order to find waves of lower intensity, we will remove the first wave with the highest intensity by subtracting it from the time series  $(y_n)$:
\[y_n^{(1)}=y_n-\frac{x_{1,max}}{1+\exp(-\frac{n-b_1}{a_1})}.
\]
Then, for the time series $(y_n^{(1)})$, we calculate its first and second differences and for the latter we perform the CWT analysis again. The above process may be repeated several times if necessary.

\begin{figure}[!h]
\centering
\begin{tabular}{ll}
\begin{subfigure}{7cm}
\centering
\includegraphics[width=7cm,height=5cm]{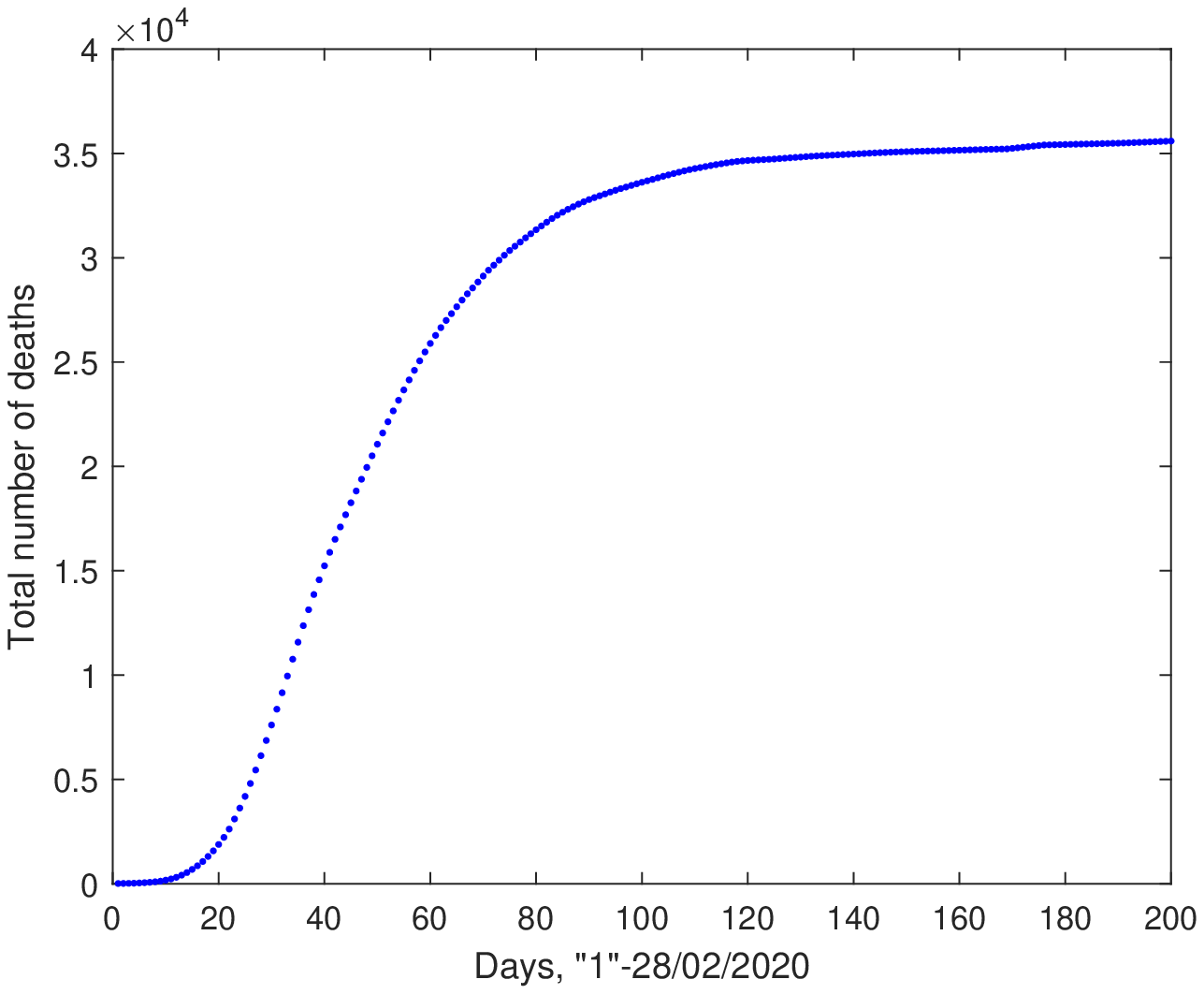}
\caption{Smoothed total number of deaths $(y_n)$}
\end{subfigure}
&
\begin{subfigure}{7cm}
\centering
\includegraphics[width=7cm,height=5cm]{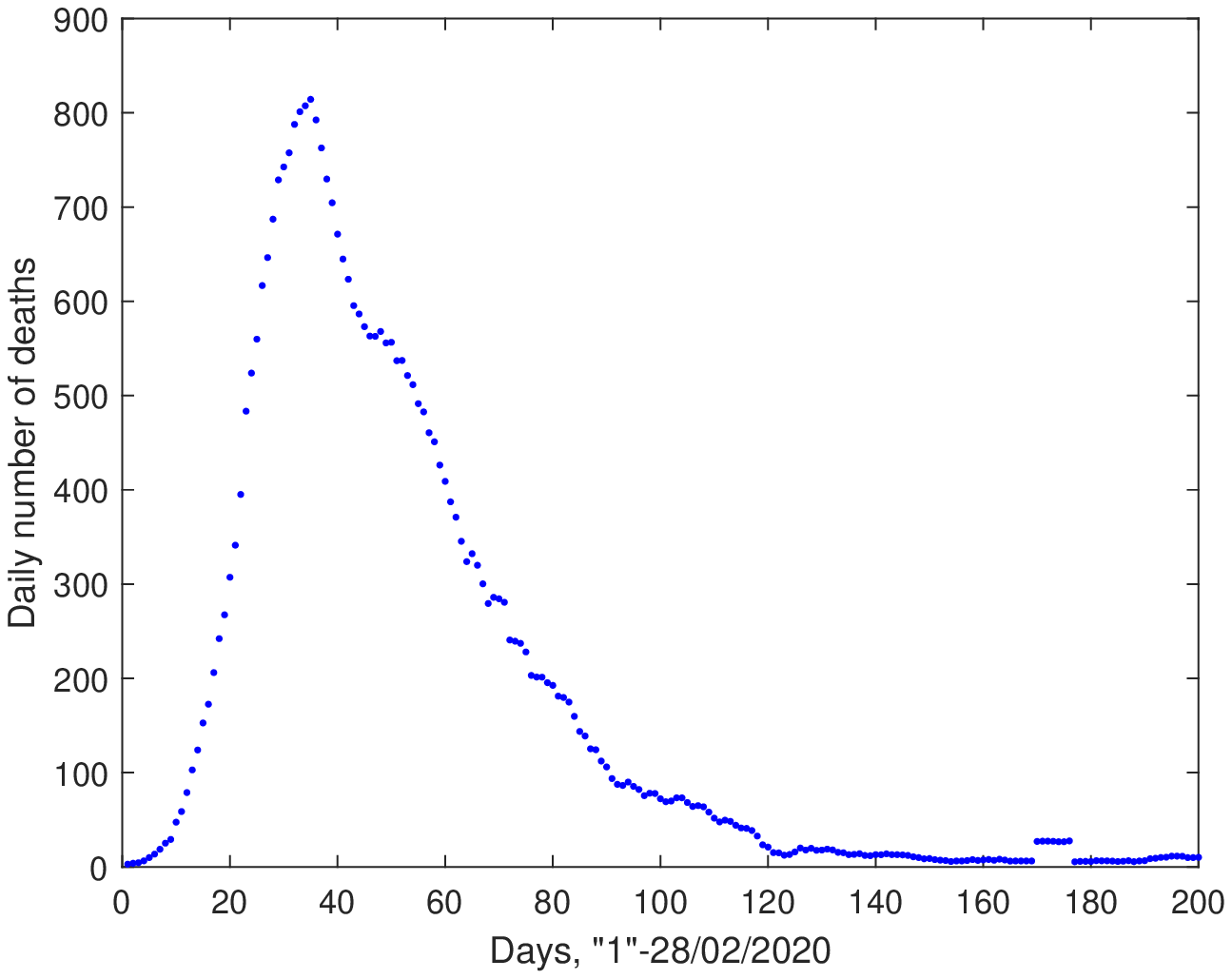}
\caption{First differences $(\Delta^{1}y_n)$}
\end{subfigure}
\end{tabular}
\centering
\begin{tabular}{ll}
\begin{subfigure}{7cm}
\centering
\includegraphics[width=7cm,height=5cm]{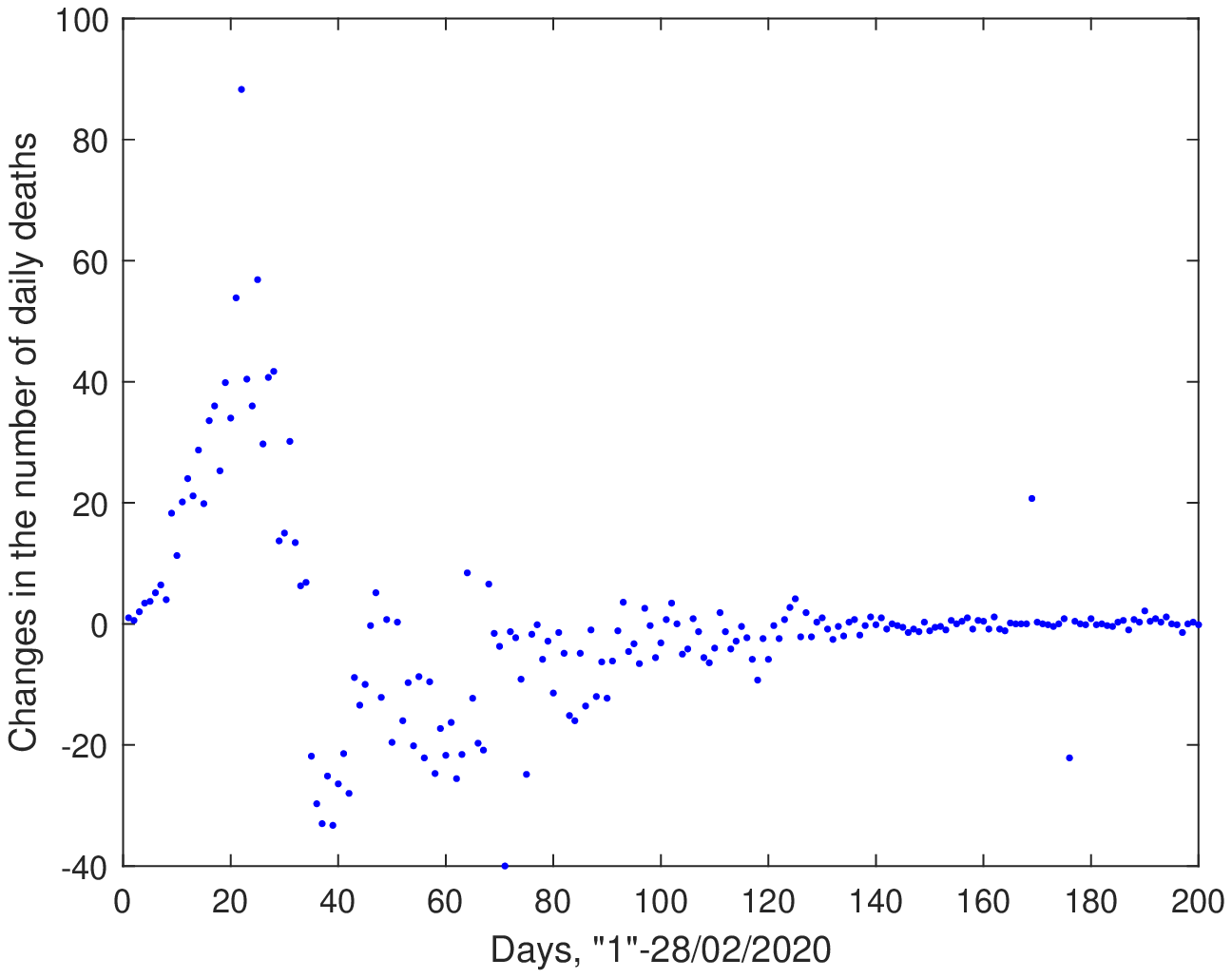}
\caption{Second differences $(\Delta^{2}y_n)$}
\end{subfigure}
&
\begin{subfigure}{7cm}
\centering
\includegraphics[width=7cm,height=5cm]{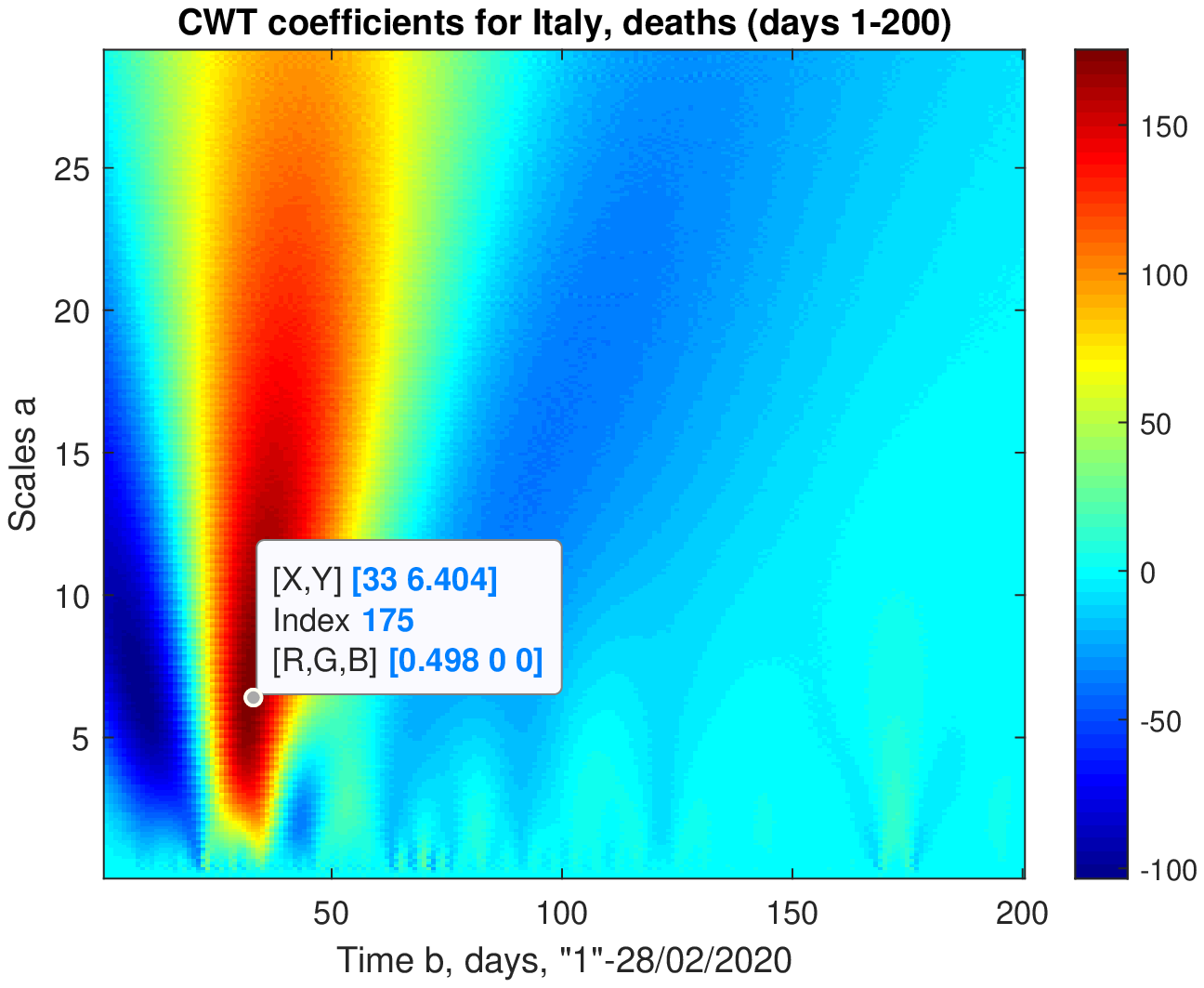}
\caption{CWT coefficients for (c)}
\end{subfigure}
\end{tabular}
\caption{`First wave'  of COVID-19 deaths in Italy, 28/02/2020--14/09/2020}
\label{fig3}
\end{figure}

Fig.~\ref{fig3}(d) shows that the logistic wave \#1 with the highest intensity has parameters $b_1=33$ and $a_1=6.4$ and the Index value is 175. The saturation level (\ref{4c}) is
\[x_{1,max}=\sqrt{30}\cdot 6.4\cdot \sqrt{6.4}\cdot 175=15,519.
\]

\begin{figure}[!h]
	\begin{center}
	 \includegraphics[height=6cm, width=8cm]{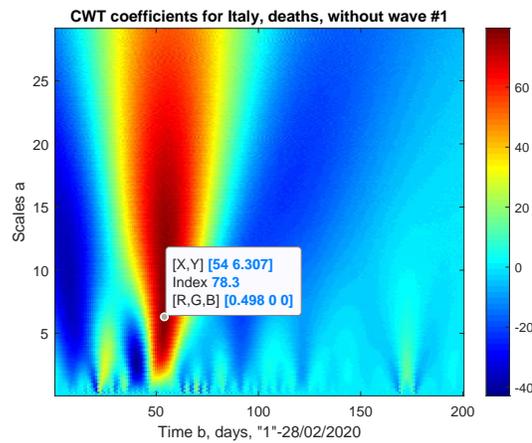}
	\end{center}
	\vspace{-7mm}
	\caption{Scalogram CWT after removing wave \#1}
	\label{fig4}
\end{figure}

After removing wave \#1, according to the procedure described above, and performing the CWT analysis, we get the scalogram shown in Fig.~\ref{fig4}, from which we read the parameters of wave \#2: $b_2=54$ and $a_2=6.3$, and calculate the saturation level (\ref{4c}) $x_{2,max}=6782$.

\begin{figure}[!h]
	\begin{center}
	 \includegraphics[height=6cm, width=8cm]{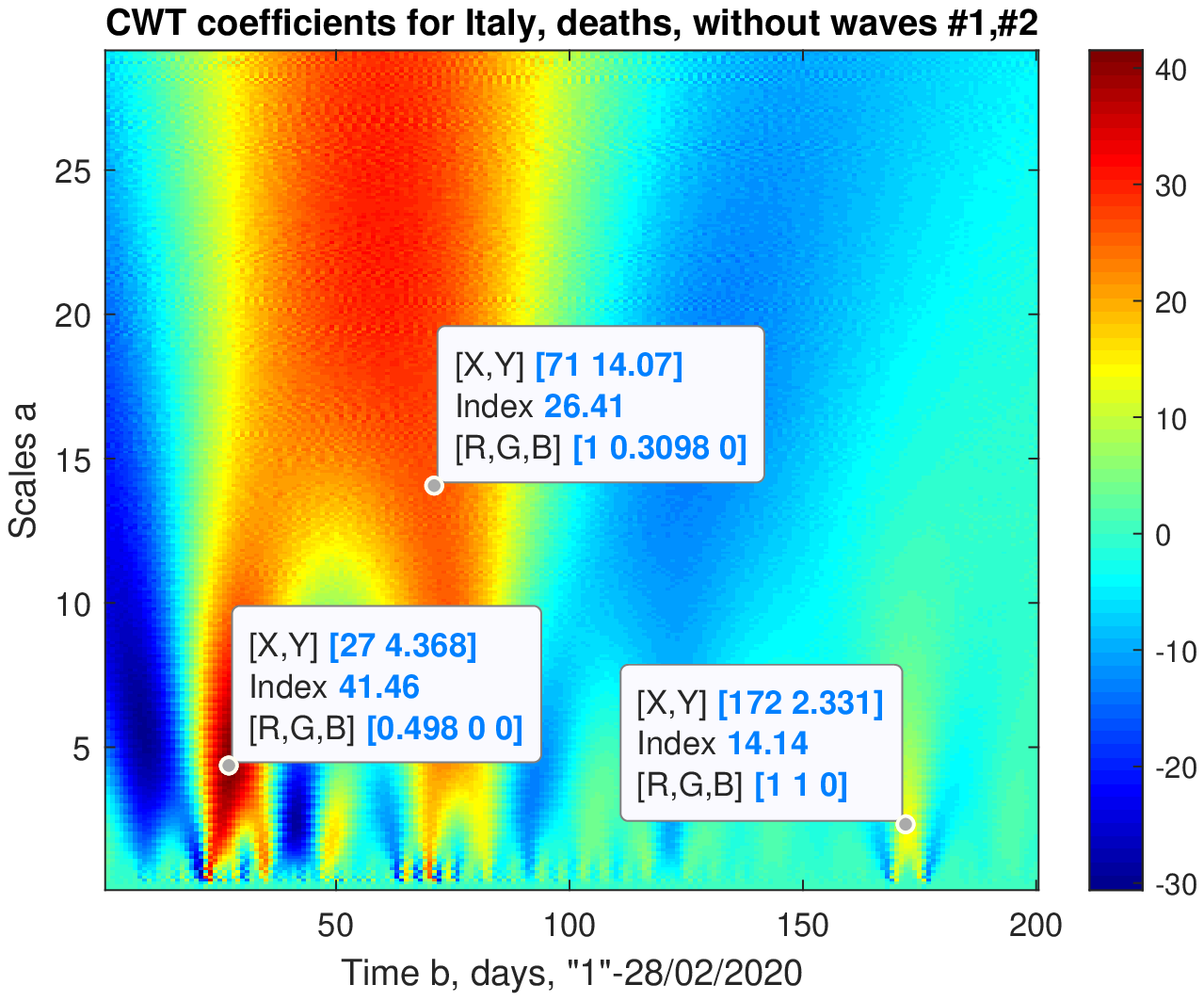}
	\end{center}
	\vspace{-7mm}
	\caption{Scalogram CWT after removing waves \#1, \#2}
	\label{fig5}
\end{figure}

In order to find logistic waves with lower intensities, we also remove wave  \#2. After performing the CWT analysis, we get the scalogram Fig.~\ref{fig5}. 

\begin{figure}[!h]
	\begin{center}
	 \includegraphics[height=6cm, width=8cm]{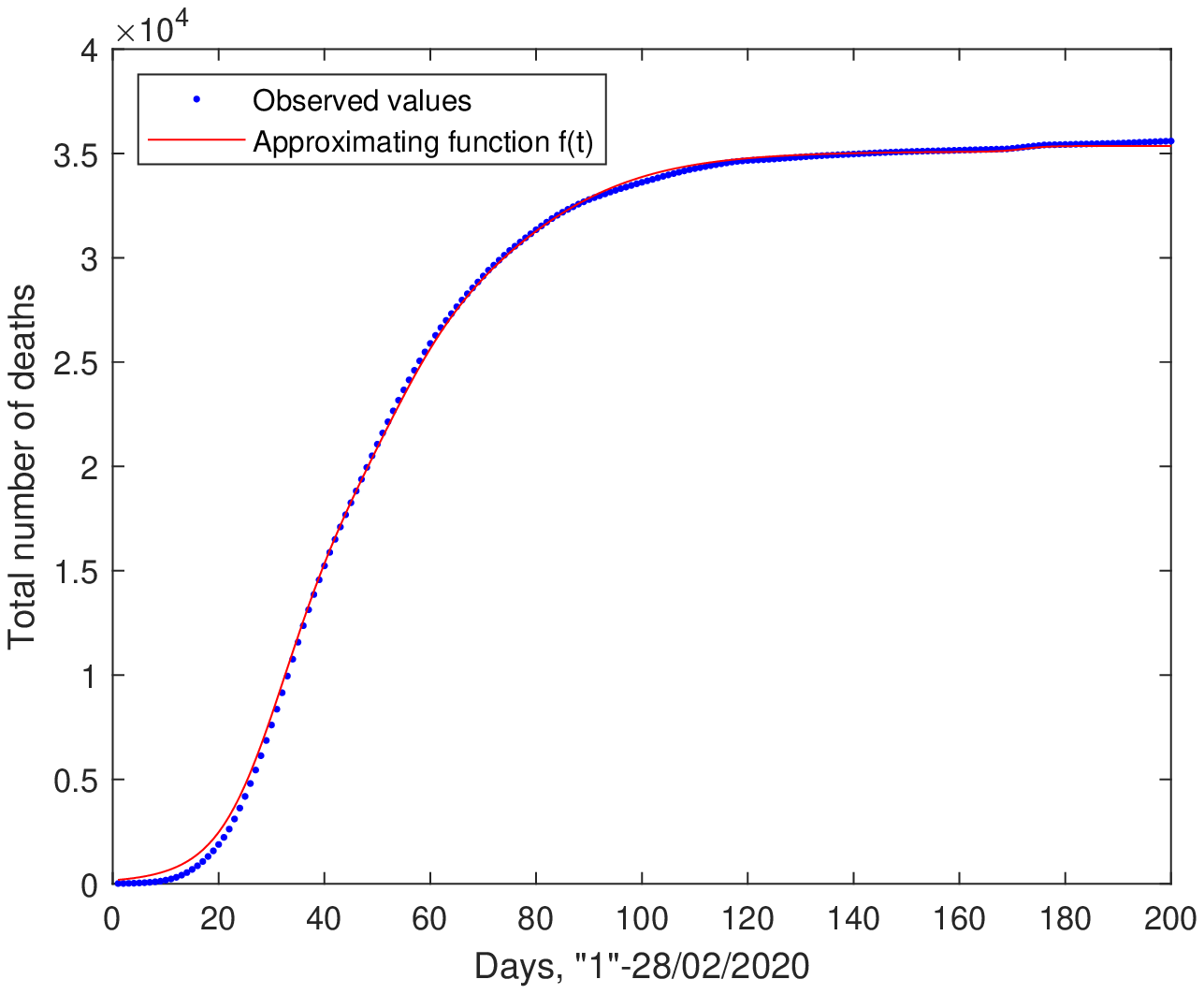}
	\end{center}
	\vspace{-7mm}
	\caption{Approximating function $f(t)$ for time series $(y_n)$}
	\label{fig6}
\end{figure}

Taking into account the three waves shown in the Fig.~\ref{fig5}, we obtain function $f(t)$ (\ref{4bb}), approximating the time series $(y_n)$, in the form of the following sum of five logistic functions (see also Fig.~\ref{fig6}):
\begin{align}\label{4d}
f(t)=\sum_{i=1}^{5}\frac{x_{i,max}}{1+\exp(-\frac{t-b_i}{a_i})}=&\frac{15,519}{1+\exp(-\frac{t-33}{6.4})}+\frac{6,782}{1+\exp(-\frac{t-54}{6.3})}
+\frac{10,692}{1+\exp(-\frac{t-71}{14.1})}\nonumber\\
&+\frac{2,098}{1+\exp(-\frac{t-27}{4.4})}+\frac{269}{1+\exp(-\frac{t-172}{2.3})}
\end{align}

The saturation levels of waves \#4 and \#5 in (\ref{4d}) were calculated according to formula (\ref{4c}). However, for wave \#3 there is no clear maximum of the Index. Therefore, and in order to compensate the influence of other, small waves, not included in the model, we have calculated its saturation level $x_{3,max}=10,692$  to minimize the value of the RMSE error:
\[\text{RMSE}=\sqrt{\frac{1}{200}\sum_{i=1}^{200}(y_n-f(n))^2}.
\]
Note that for model (\ref{4d}) $\text{RMSE}=231.02  $

\section{Conclusions}\label{sec4}
In this paper we deal with logistic wavelets and their normalization in $L^{2}(R)$ space. We then use them to study the spread of the first wave of COVID-19 deaths in Italy in 2020. It turned out that this wave, although asymmetric, can be described by the sum of five logistic functions (curves). Wave \#3 lasted as long as (\ref{4cc}), $7.32\cdot 14.1=103$ days, while waves \#1 and \#2 were of similar length, 47 and 46 days respectively. The peaks of daily deaths for waves \#1 and \#2 were $b_2-b_1=54-33= 21$ days apart. So, against the background of the long wave \#3, occurred waves \#1 and \#2, about half as long. Wave \#4 arrived a few days earlier than wave \#1, but was much less intense. Wave \#5 was a single pulse, with a low value of the saturation level.

{\large\textbf{Funding statement}}

The research of the author was partially funded by the 'IDUB against COVID-19' project granted by the Warsaw University of Technology (Warsaw, Poland) under the program Excellence Initiative: Research University (IDUB), grant no 1820/54/201/2020.

\vspace{7mm}

{\large\textbf{Conflict of Interests}}  

The author declares that  there is no any conflict of interest in the submitted manuscript.

\end{document}